\documentclass[11pt]{amsart}
\usepackage{amssymb}
\usepackage[all]{xy}
\usepackage[usenames]{color}

\title[Computing Hilbert-Siegel modular forms]{Computing genus $2$ Hilbert-Siegel modular forms over $\Q(\sqrt{5})$ via the Jacquet-Langlands correspondence}

\author{Clifton Cunningham}
\author{Lassina Demb\'el\'e}

\date{\today}

\address{Department of Mathematics, University of Calgary}
\email{cunning@math.ucalgary.ca} 
\address{Institut f\"{u}r Experimentelle Mathematik, Universit\"{a}t Duisburg-Essen}
\email{lassina.dembele@uni-duisburg-essen.de}

\dedicatory{}

\newtheorem{thm}{Theorem}
\newtheorem{prop}[thm]{Proposition}

\newtheorem{rem}{Remark}

\newtheorem{lem}{Lemma}

\newtheorem{conj}[thm]{Conjecture}

\theoremstyle{definition}

\theoremstyle{remark}

\numberwithin{equation}{section}

\newenvironment{prf}{\begin{proof}[\bf{Proof}]}{\end{proof}}

\newcommand{\field}[1]{\mathbb{#1}}
\newcommand{\Q}{\field{Q}}
\newcommand{\C}{\field{C}}
\newcommand{\R}{\field{R}}

\newcommand{\Z}{\field{Z}}

\newcommand{\F}{\field{F}}

\newcommand{\ceq}{{\, :=\, }}
\newcommand{\tq}{{\ \big\vert\ }}

\begin{document}

\begin{abstract} 
In this paper we present an algorithm for computing Hecke eigensystems of Hilbert-Siegel cusp forms over real quadratic fields of narrow class number one. We give some illustrative examples using the quadratic field $\Q(\sqrt{5})$. In those examples, we identify Hilbert-Siegel eigenforms that are possible lifts from Hilbert eigenforms.
\end{abstract}

\maketitle


\section*{\bf Introduction}

Let $F$ be a real quadratic field of narrow class number one and let $B$ be the unique (up to isomorphism) quaternion algebra over $F$ which is ramified at both archimedean places of $F$ and unramified everywhere else. Let $\mathbf{GU}_2(B)$ be the unitary similitude group of $B^{\oplus 2}$. This is the set of $\Q$-rational points of an algebraic group $G^B$ defined over $\Q$. The group $G^B$ is an inner form of $G :=\mathrm{Res}_{F/\Q}(\mathbf{GSp}_4)$  such that $G^B(\R)$ is compact modulo its centre. {(These notions are reviewed at the beginning of Section~\ref{section: automorphicforms}.)}

In this paper we develop an algorithm which computes automorphic forms on $G^B$ in the following sense: given an ideal $N$ in $\mathcal{O}_F$ and an integer $k$ greater than {$2$}, the algorithm returns the Hecke eigensystems of all automorphic forms $f$ of level $N$ and parallel weight $k$. More precisely, given a prime $\mathfrak{p}$ in $\mathcal{O}_F$, the algorithm returns the Hecke eigenvalues of $f$ at $\mathfrak{p}$, and hence the Euler factor $L_\mathfrak{p}(f,s)$, for each eigenform $f$ of level $N$ and parallel weight $k$. The algorithm is a generalization of the one developed in \cite{dembele1} to the genus 2 case. Although we have only described the algorithm in the case of a real quadratic field in this paper, it should be clear from our presentation that it can be adapted to any totally real number field of narrow class number one. 

{The Jacquet-Langlands Correspondence of the title refers to the conjectural} map $JL:\,\Pi(G^B)\to \Pi(G)$ from automorphic representations of $G^B$ to automorphic representations of $G$, which is injective, matches L-functions and enjoys other properties compatible with the principle of functoriality; in particular, the image of the Jacquet-Langlands Correspondence is to be contained in the space of holomorphic automorphic representations. If we admit this conjecture, then the algorithm above provides a way to produce examples of cuspidal Hilbert-Siegel modular forms of genus 2 over $F$ and allows us to compute the L-factors of the corresponding automorphic representations for arbitrary finite primes $\mathfrak{p}$ of $F$. 

In fact, we are also able to use these calculations to provide evidence for the Jacquet-Langlands Correspondence itself by comparing the Euler factors we find with those of known Hilbert-Siegel modular forms obtained by lifting. This we do in the final section of the paper where we observe that some of the Euler factors we compute match those of lifts of Hilbert modular forms, for the primes we computed. Although this does not definitively establish that these Hilbert-Siegel modular forms are indeed lifts, in principle  one can establish equality in this way, using an analogue of the Sturm bound.

The first systematic approach to Siegel modular forms from a computational viewpoint is due to Skoruppa \cite{skoruppa1} who used Jacobi symbols to generate spaces of such forms. His algorithm, which has been extensively exploited by Ryan \cite{ryan}, applies only to the case of full level structure. More recently, Faber and van der Geer \cite{FvdG1} and \cite{FvdG2} also produced examples of Siegel modular forms by counting points on hyperelliptic curves of genus 2; again their results are available only in the full level structure case. The most substantial progress toward the computation of Siegel modular forms for proper level structure is by Gunnells \cite{gunnells} who extended the theory of modular symbols to the symplectic group $\mathbf{Sp}_4/\Q$. However, this work does not see the cuspidal cohomology, which is the only part of the cohomology which is relevant to arithmetic geometric applications.
To the best of our knowledge, there are no numerical examples of Hilbert-Siegel modular forms for proper level structure in the literature, with the exception of those produced from liftings of Hilbert modular forms.

The outline of the paper is as follows. In Section~\ref{section: automorphicforms} we recall the basic properties of Hilbert-Siegel modular forms and algebraic automorphic forms together with the Jacquet-Langlands Correspondence. In Section~\ref{section: algorithm} we give a detailed description of our algorithm. Finally, in Section~\ref{section: examples} we present numerical results for the quadratic field $\Q(\sqrt{5})$.

\medskip\noindent{\bf Acknowledgements.} {During the course of the preparation of this paper, the second author had helpful email exchanges with several people including Alexandru Ghitza, David Helm, Marc-Hubert Nicole, David Pollack, Jacques Tilouine and Eric Urban. The authors wish to thank them all. Also, we would like to thank William Stein for allowing us to use the SAGE computer cluster at the University of Washington. And finally, the second author would like to thank the PIMS institute for their postdoctoral fellowship support, and the University of Calgary for {its} hospitality.}


\section{\bf Hilbert-Siegel modular forms and the Jacquet-Langlands correspondence}\label{section: automorphicforms}

Throughout this paper, $F$ denotes a real quadratic field of narrow class number one. The two archimedean places of $F$ and the real embeddings of $F$ will both be denoted $v_0$ and $v_1$. For every $a\in F$, we write $a_0$ (resp. $a_1$) for the image of $a$ under $v_0$ (resp. $v_1$). The ring of integers of $F$ is denoted by $\mathcal{O}_F$. For every prime ideal $\mathfrak{p}$ in $\mathcal{O}_F$, the completion of $F$ and $\mathcal{O}_F$ at $\mathfrak{p}$ will be denoted by $F_\mathfrak{p}$ and $\mathcal{O}_{F_\mathfrak{p}}$, respectively. 

Let $B$ be the unique (up to isomorphism) totally definite quaternion algebra over $F$ which is unramified at all finite primes of $F$.  We fix a maximal order $\mathcal{O}_B$ of $B$. Also, we choose a splitting field $K/F$ of B that is Galois over $\Q$ and such that there exists an isomorphism $j:\,\mathcal{O}_B\otimes_\Z\mathcal{O}_K \cong \mathbf{M}_2(\mathcal{O}_K)\oplus \mathbf{M}_2(\mathcal{O}_K)$, where $\mathbf{M}_2(A)$ denotes the ring of $2\times 2$-matrices with entries from a ring $A$. For every finite prime $\mathfrak{p}$ in $F$, we fix an isomorphism $B_\mathfrak{p}\cong \mathbf{M}_2(F_\mathfrak{p})$ which restricts to an isomorphism from $\mathcal{O}_{B,\,\mathfrak{p}}$ onto $\mathbf{M}_2(\mathcal{O}_{F_\mathfrak{p}})$.


The algebraic group $G=\mathrm{Res}_{F/\Q}(\mathbf{GSp}_{4})$ is defined as follows. For any $\Q$-algebra $A$, the set of $A$-rational points {of $G$} is given by
\[
G(A)= \left\{\gamma \in \mathbf{GL}_{4}(A\otimes_{\Q}F) \Big\vert \begin{array}{c}\gamma J_{2}\gamma^t = \nu_{G}(\gamma) J_{2} \\ \nu_G(\gamma)\in (A\otimes_\Q F)^\times \end{array}\right\},
\]
where 
\[
J_{2}= \begin{pmatrix} 0 & \mathbf{1}_2 \\ -\mathbf{1}_2 & 0 \end{pmatrix}.
\]
This group admits an integral model {with} $A$-rational points for every $\Z$-algebra $A$ given by
\[
 G_{\Z}(A)= \left\{\gamma \in \mathbf{GL}_{4}(A\otimes_{\Z}\mathcal{O}_F)\Big\vert \begin{array}{c}\gamma J_{2}\gamma^t = \nu_{G}(\gamma) J_{2} \\ \nu_G(\gamma)\in (A\otimes_\Z \mathcal{O}_F)^\times\end{array} \right\}.
\]
For any $\Q$-algebra $A$, the conjugation on $B$ extends in a natural way to the matrix algebra $\mathbf{M}_2(B\otimes_\Q A)$. 

The algebraic group $G^B/\Q$ is defined {as follows. For any $\Q$-algebra $A$, the set of $A$-rational points of $G^B$ is given by}
\[
G^B(A)= \left\{\gamma \in \mathbf{M}_{2}(B\otimes_{\Q}A)\Big\vert \begin{array}{c} \gamma\bar{\gamma}^t = \nu_{G^B}(\gamma)\mathbf{1}_2 \\ \nu_{G^B}(\gamma)\in (A\otimes_\Q F)^\times\end{array}\right\}.
\]
This group also admits an integral model with $A$-rational points for every $\Z$-algebra given by
\[
G_{\Z}^B(A)= \left\{\gamma \in \mathbf{M}_2(\mathcal{O}_B\otimes_{\Z}A) \Big\vert \begin{array}{c} \gamma\bar{\gamma}^t= \nu_{G^B}(\gamma)\mathbf{1}_2 \\ \nu_{G^B}(\gamma)\in (A\otimes_\Z \mathcal{O}_F)^\times\end{array}\right\}.
\]

The group $G^B/\Q$ is an inner form of $G/\Q$ such that $G^B(\R)$ is compact modulo its center. 
Combining the isomorphism $j$ (see above) with conjugation by a permutation matrix, we obtain an isomorphism $G^B_\Z(\mathcal{O}_K)\cong G_\Z(\mathcal{O}_K)$, which we fix from now on.
For every prime ideal $\mathfrak{p}$ in $F$, the splitting of $G^B$ at $\mathfrak{p}$ amounts to the splitting of the quaternion algebra $B$ at $\mathfrak{p}$; we refer to \cite{dembele1} for further details.

By the choice of the quaternion algebra $B$, we have $G^B(\hat{\Q})\cong G(\hat{\Q})$. (We denote the finite ad\`eles of $\Q$ (resp. $\Z$) by $\hat{\Q}$ (resp. $\hat{\Z}$)).

\subsection{\bf Hilbert-Siegel modular forms}
We fix an integer $k\ge 3$ and, for simplicity, we restrict ourselves to Hilbert-Siegel modular forms of parallel weight $k$. The real embeddings $v_0$ and $v_1$ of $F$ extend to $G(\Q)=\mathbf{GSp}_4(F)$ in a natural way.
We denote by $\mathbf{GSp}_4^+(F)$ the subgroup of elements $\gamma$ with totally positive similitude factor $\nu_G(\gamma)$. We recall that the Siegel upper-half plane of genus 2 is defined by
$$
\mathfrak{H}_2=\{\gamma\in\mathbf{GL}_2(\C) \tq \gamma^t=\gamma\,\mbox{and}\,\mathrm{Im}(\gamma)\, {\mbox{is positive definite}\,} \}.
$$
{We also recall} that $\mathbf{GSp}_4^{+}(F)$ acts on $\mathfrak{H}_2^2$ by
\[
\begin{pmatrix}a&b\\  c& d\end{pmatrix} (\tau_0,\,\tau_1) \ceq\left((a_0\tau_0+b_0)(c_0\tau_0+d_0)^{-1},\, (a_1\tau_1+b_1)(c_1\tau_1+d_1)^{-1}\right).
\]
This induces an action on the space of functions $f:\mathfrak{H}_2^2\to\C$ by 
$$
\forall\gamma=\begin{pmatrix}a&b\\  c& d\end{pmatrix}, \qquad f|_k\gamma(\tau)= \prod_{i=0}^1\frac{\nu_G(\gamma_i)^{k/2}}{\det(c_i\tau_i+d_i)^k} f(\tau).
$$

Let $N$ be an ideal in $\mathcal{O}_F$ and set 
$$
\Gamma_0(N)=\left\{\begin{pmatrix}a & b \\ c & d \end{pmatrix} \in \mathbf{GSp}_{4}^{+}(\mathcal{O}_F) \tq c\equiv 0 (N) \right\}.
$$
A {\bf Hilbert-Siegel modular form} of level $N$ and parallel weight $k$ is a holomorphic function $f:\,\mathfrak{H}_2^2\to\C$ such that
$$
\forall \gamma\in\Gamma_0(N), \qquad f|_k\gamma=f.
$$
The space of Hilbert-Siegel modular forms of parallel weight $k$ and level $N$ is denoted $M_k(N)$. Each $f\in M_k(N)$ admits a Fourier expansion, which by the Koecher principle takes the form
$$
\forall \tau\in\mathfrak{H}_2^2, \qquad f(\tau)=\sum_{\{Q\}\cup \{0\}}a_Q e^{2\pi i\mathrm{Tr}(Q\tau)},
$$ 
where $Q\in  \mathbf{M}_2(F)$ runs over all symmetric totally positive and semi-definite matrices. A Hilbert-Siegel modular forms $f$ is a {\bf cusp form} if, for all $\gamma\in\mathbf{GSp}_4^{+}(F)$, the {constant term} in the Fourier expansion of $f|_k\gamma$ is zero. The space of Hilbert-Siegel cusp forms is denoted $S_k(N)$.

\subsection{\bf The Hecke algebra} The space $S_k(N)$ comes equipped with a Hecke action, which we now recall. Take $u\in\mathbf{GSp}_4^{+}(F)\cap \mathbf{M}_4(\mathcal{O}_F)$, and write the finite disjoint union
$$
\Gamma_0(N)u\Gamma_0(N)=\coprod_{i}\Gamma_0(N)u_i.
$$ 
Then the Hecke operator $[\Gamma_0(N)u\Gamma_0(N)]$ on $S_k(N)$ is given by
$$
[\Gamma_0(N)u\Gamma_0(N)]f=\sum_{i}f|_ku_i.
$$
Let $\mathfrak{p}$ be a prime ideal in $\mathcal{O}_F$ and let $\pi_\mathfrak{p}$ be a totally positive generator of $\mathfrak{p}$; let $T_1(\mathfrak{p})$ and $T_2(\mathfrak{p})$ be the Hecke operators corresponding to the double $\Gamma_0(N)$-cosets of the {symplectic similitude} matrices
\[
\begin{pmatrix}
1 & 0 & 0 & 0 \\ 
0 & 1 & 0 & 0 \\
0 & 0 & \pi_{\mathfrak{p}} & 0 \\
0 & 0 & 0 & \pi_{\mathfrak{p}} 
\end{pmatrix}
\quad\text{and}\quad 
\begin{pmatrix}
1 & 0 & 0 & 0 \\ 
0 & \pi_{\mathfrak{p}} & 0 & 0 \\
0 & 0 & \pi_{\mathfrak{p}}^2 & 0 \\
0 & 0 & 0 & \pi_{\mathfrak{p}} 
\end{pmatrix},
\]
respectively. {(We remind the reader of the symplectic form $J_2$ fixed at the beginning of Section~\ref{section: automorphicforms}.)} The Hecke algebra $\mathbf{T}_k(N)$ is the $\Z$-algebra generated by the operators $T_1(\mathfrak{p})$ and $T_2(\mathfrak{p})$, where $\mathfrak{p}$ runs over all primes not dividing $N$.

\subsection{\bf Algebraic Hilbert-Siegel autormorphic forms} We only consider level structure of Siegel type. Namely, we define the compact open subgroup $U_0(N)$ of
$G(\hat{\Q})$ by
$$
U_0(N)=\prod_{\mathfrak{p}\nmid N}\mathbf{GSp}_4(\mathcal{O}_{F_\mathfrak{p}})\times\prod_{\mathfrak{p}\mid N}U_0(\mathfrak{p}^{e_\mathfrak{p}}),
$$ 
where
$N=\prod_{\mathfrak{p}\mid N}\mathfrak{p}^{e_\mathfrak{p}}$ and
$$
U_0(\mathfrak{p}^{e_\mathfrak{p}}):=\left\{\begin{pmatrix}a&b\\ c&d \end{pmatrix}\in\mathbf{GSp}_4(\mathcal{O}_{F_\mathfrak{p}}) \tq c\equiv 0\mod\mathfrak{p}^{e_\mathfrak{p}}\right\}.
$$ 

The {\bf weight representation} is defined as follows. Let $L_{k}$ be the representation of $\mathbf{GSp}_4(\C)$ of highest weight $(k-3,\,k-3)$. We let $V_k=L_k\otimes L_k$ and define the complex representation
$(\rho_k,\,V_k)$ by
$$
\rho_k:\,G^B(\R)\longrightarrow\mathbf{GL}(V_k),
$$
where the action on the first factor is via $v_0$, and the action on the second one is via $v_1$.

The space of {\bf algebraic Hilbert-Siegel modular forms} of weight $k$ and level $N$ is given by
$$
M_k^B(N):=\left\{f:\,G^B(\hat{\Q})/U_0(N)\to V_k \tq \forall \gamma\in G^B(\Q), f|_k\gamma=f \right\},
$$
where $f|_k\gamma(x)=f(\gamma x)\gamma$, for all $x\in G^B(\hat{\Q})/U_0(N)$. When $k=3$, we let 
$$
I_k^B(N):=\left\{f:\,G^B(\Q)\backslash G^B(\hat{\Q})/U_0(N)\to\C \tq \text{$f$ is constant} \right\}.
$$ 
Then, the space of {\bf algebraic Hilbert-Siegel cusp forms} of weight $k$ and level $N$ is defined by
\begin{eqnarray*}
S_k^B(N):=\left\{\begin{array}{lll} 
M_k^B(N)&\mbox{if}& k>3,\\
&& \\
M_k^B(N)/I_k^B(N)&\mbox{if}&k=3.\\
\end{array}\right.
\end{eqnarray*}

The action of the Hecke algebra on $S_{k}^B(N)$ is given as follows. For any $u\in G(\hat{\Q})$, write the
finite disjoint union 
$$
U_0(N)uU_0(N)=\coprod_{i}u_iU_0(N),
$$ 
and define
\begin{eqnarray*}
[U_0(N)uU_0(N)]:\,S_{k}^B(N)&\to&S_{k}^B(N)\\
f&\mapsto&f|_{k}[U_0(N)uU_0(N)],
\end{eqnarray*}
by
\[
f|_{k}[U_0(N)uU_0(N)](x)=\sum_{i}f(xu_i),\,x\in G(\hat{\Q}).
\]
For any prime $\mathfrak{p}\nmid N$, let $\varpi_\mathfrak{p}$ be a local uniformizer at $\mathfrak{p}$. The local Hecke algebra at $\mathfrak{p}$ is generated by the Hecke operators $T_1(\mathfrak{p})$ and $T_2(\mathfrak{p})$ corresponding to the double $U_0(N)$-cosets $\Delta_1(\mathfrak{p})$ and $\Delta_2(\mathfrak{p})$ of the matrices 
\[
\begin{pmatrix}
1 & 0 & 0 & 0 \\ 
0 & 1 & 0 & 0 \\
0 & 0 & \varpi_{\mathfrak{p}} & 0 \\
0 & 0 & 0 & \varpi_{\mathfrak{p}} 
\end{pmatrix}
\quad\text{and}\quad 
\begin{pmatrix}
1 & 0 & 0 & 0 \\ 
0 & \varpi_{\mathfrak{p}} & 0 & 0 \\
0 & 0 & \varpi_{\mathfrak{p}}^2 & 0 \\
0 & 0 & 0 & \varpi_{\mathfrak{p}} 
\end{pmatrix},
\]
respectively. We let $\mathbf{T}_{k}^B(N)$ be the Hecke algebra generated by $T_1(\mathfrak{p})$ and $T_2(\mathfrak{p})$ for all primes $\mathfrak{p}\nmid N$.

\subsection{\bf The Jacquet-Langlands Correspondence}\label{subsection: JL}
 The Hecke modules $S_k(N)$ and $S_k^B(N)$ are related by the following conjecture known as the Jacquet-Langlands Correspondence {for symplectic similitude groups}.

\begin{conj}\label{jacquetlanglands} The Hecke algebras $\mathbf{T}_k(N)$ and $\mathbf{T}_k^B(N)$ are isomorphic and there is a compatible isomorphism of Hecke modules
$$S_k(N)\stackrel{\sim}{\longrightarrow} S_k^B(N).$$
\end{conj}

{It is common, but perhaps not entirely accurate, to attribute this conjecture to Jacquet-Langlands. To the best of our knowledge, the correspondence in this form was first discussed by Ihara \cite{ihara} in the case $F=\Q$.
In \cite{ibukiyama}, Ibukiyama provided some numerical evidence. 
On the other hand, it is appropriate to refer to Conjecture~\ref{jacquetlanglands} as the Jacquet-Langlands Correspondence (for $GSp(4)$) since it is an analogue of the Jacquet-Langlands Correspondence (for $GL(2)$) which relates automorphic representations of the multiplicative group of a quaternion algebra with certain automorphic representations of $GL(2)$ (see \cite{JL}). Both correspondences are, in turn, special consequences of the principle of functoriality, as expounded by Langlands. 
Finally, it appears that Conjecture~\ref{jacquetlanglands} may soon be a theorem due to the work of \cite{sorenson} and the forthcoming book by James Arthur on automorphic representations of classical groups.}

\section{\bf The Algorithm}\label{section: algorithm}

In this section, we present the algorithm we used in order to compute the Hecke module of (algebraic) Hilbert-Siegel modular forms. The main assumption in this section is that the class number of the principal genus of $G^B$ is 1. ({We refer to \cite{dembele3} to see how one  can relax this condition on the class number.}) We recall that since $B$ is totally definite,  $G^B$ satisfies Proposition 1.4 in Gross \cite{gross1}. Thus the group $G^B(\R)$ is compact modulo its centre, and $\Gamma=G^B(\Z)/\mathcal{O}_F^\times$ is {\it finite}.

For any prime $\mathfrak{p}$ in $F$, let $\F_\mathfrak{p}=\mathcal{O}_F/\mathfrak{p}$ be the residue field at $\mathfrak{p}$ and define the reduction map
\begin{eqnarray*}
\mathbf{M}_2(\mathcal{O}_{B,\,\mathfrak{p}})&\to&\mathbf{M}_4(\F_\mathfrak{p})\\
g&\mapsto&\tilde{g},
\end{eqnarray*} 
{where we use the splitting of $\mathcal{O}_{B,\mathfrak{p}}$ that was fixed at the beginning of Section~\ref{section: automorphicforms}.}
Now, choose a totally positive generator $\pi_\mathfrak{p}$ of $\mathfrak{p}$ and put
\[
\begin{aligned}
\Theta_1(\mathfrak{p}) 
&:=\Gamma\backslash\left\{u\in \mathbf{M}_2(\mathcal{O}_B) \tq u\bar{u}^t=\pi_{\mathfrak{p}}\mathbf{1}_2 \text{and } \mathrm{\rm rank}(\tilde{g})=2\right\},\\
\Theta_2(\mathfrak{p})
&:=\Gamma\backslash\left\{u\in \mathbf{M}_2(\mathcal{O}_B) \tq
u\bar{u}^t=\pi_{\mathfrak{p}}^2\mathbf{1}_2\,
\text{and }\,\mathrm{\rm rank}(\tilde{g})=1\right\}.
\end{aligned}
\]
We let $\mathcal{H}_0^2(N)=G(\hat{\Z})/U_0(N)$.  Then the group $\Gamma$ acts on $\mathcal{H}_0^2(N)$, thus on the space of functions $f:\,\mathcal{H}_0^2(N)\to V_k$ by
$$
\forall x\in\mathcal{H}_0^2(N), \forall \gamma\in\Gamma, \qquad
f|_k\gamma(x):=f(\gamma x)\gamma.
$$ 

\begin{thm}\label{thm1} There is an isomorphism of Hecke modules
$$
M_{k}^B(N)\stackrel{\sim}{\longrightarrow}\left\{f:\,\mathcal{H}_0^2(N)\to V_{k} \tq f|_{k}\gamma=f,\,
\gamma\in\Gamma\right\},
$$ 
where the Hecke action on the right hand side is given by
$$
\begin{aligned}
f|_{k}T_1(\mathfrak{p}) & =\sum_{u\in\Theta_1(\mathfrak{p})}f|_{k}u, \\
f|_{k}T_2(\mathfrak{p}) & =\sum_{u\in\Theta_2(\mathfrak{p})}f|_{k}u.
\end{aligned}
$$
\end{thm}

\begin{prf} 
The canonical map
$$\phi:\,G^B(\Z)\backslash G^B(\hat{\Z})/U_0(N)\rightarrow G^B(\Q)\backslash G^B(\hat{\Q})/U_0(N)$$
is an injection. Making use of the fact that the class number in the principal genus of $G^B$ is one ($G^B(\hat{\Q})=G^B(\Q)G_{\Z}^B(\hat{\Z})$), we see that $\phi$ is in fact a bijection. Since each element $f\in M_{k}^B(N)$ is determined by its values on a set of coset representatives of $G^B(\Q)\backslash G^B(\hat{\Q})/U_0(N)$, the map $\phi$ induces an isomorphism of complex vector spaces
\begin{eqnarray*}
M_{k}^B(N)&\stackrel{\sim}{\longrightarrow}&\left\{f:\,\mathcal{H}_0^2(N)\to V_{k} \tq f|_{k}\gamma=f,\, \gamma\in\Gamma\right\}\\
f&\longmapsto&f\circ\phi.
\end{eqnarray*}
We make this into a Hecke module isomorphism by defining the Hecke action on the right hand side as indicated in the statement of the theorem.
\end{prf}

In the rest of this section, we explain the main steps of the algorithm provided by Theorem~\ref{thm1}.

\subsection{\bf The quotient $\mathcal{H}_0^2(N)$}\label{quotient} 
Keeping the notations of the previous section, we recall that $N=\prod_{\mathfrak{p}\mid N}\mathfrak{p}^{e_\mathfrak{p}}$. Let $\mathfrak{p}$ be a prime dividing $N$
and consider the rank $4$ free $\left(\mathcal{O}_{F_\mathfrak{p}}/\mathfrak{p}^{e_\mathfrak{p}}\right)$-module
$L=\left(\mathcal{O}_{F_\mathfrak{p}}/\mathfrak{p}^{e_\mathfrak{p}}\right)^4$ endowed with the symplectic pairing $\langle\,,\,\rangle$
given by the matrix
$$J_2=\begin{pmatrix}0&\mathbf{1}_2\\ -\mathbf{1}_2&0 \end{pmatrix},$$ where $\mathbf{1}_2$ is the identity matrix in
$\mathbf{M}_2(\mathcal{O}_{F_\mathfrak{p}}/\mathfrak{p}^{e_\mathfrak{p}})$. Let $M$ be a rank 2
$\left(\mathcal{O}_{F_\mathfrak{p}}/\mathfrak{p}^{e_\mathfrak{p}}\right)$-submodule which is a direct factor in $L$. We say that $M$ is
{isotropic} if $\langle u,\,v\rangle =0$ for all $u,\,v\in M$. We recall that
$\mathbf{GSp}_4(\mathcal{O}_{F_\mathfrak{p}})$ acts transitively on the set of {rank 2, isotropic
$\left(\mathcal{O}_{F_\mathfrak{p}}/\mathfrak{p}^{e_\mathfrak{p}}\right)$-submodules of $L$} and that the stabilizer of the submodule
generated by $e_1=(1,\,0,\,0,\,0)^{\tiny T}$ and $e_2=(0,\,1,\,0,\,0)^{\tiny T}$ is $U_0(\mathfrak{p}^{e_\mathfrak{p}})$.
The quotient
$\mathcal{H}_0^2(\mathfrak{p}^{e_\mathfrak{p}})=\mathbf{GSp}_4(\mathcal{O}_{F_\mathfrak{p}})/U_0(\mathfrak{p}^{e_\mathfrak{p}})$ is {the set of rank 2, isotropic $\left(\mathcal{O}_{F_\mathfrak{p}}/\mathfrak{p}^{e_\mathfrak{p}}\right)$-submodules of $L$.}
Via the reduction map $\hat{\mathcal{O}}_F\to\mathcal{O}_F/N$, the quotient
$G_{\Z}(\hat{\Z})/U_0(N)$ can be identified with the product $$\mathcal{H}_0^2(N)=\prod_{\mathfrak{p}\mid N}\mathcal{H}_0^2(\mathfrak{p}^{e_\mathfrak{p}}).$$
The cardinality of $\mathcal{H}_0^2(N)$ is extremely useful and is determined using the following lemma.

\begin{lem}Let $\mathfrak{p}$ be a prime in $F$ and $e_\mathfrak{p}\ge 1$ an integer. Then, the cardinality of the set $\mathcal{H}_0^2(\mathfrak{p}^{e_\mathfrak{p}})$ is given by
$$\#\mathcal{H}_0^2(\mathfrak{p}^{e_\mathfrak{p}})=\mathbf{N}(\mathfrak{p})^{3(e_\mathfrak{p}-1)}(\mathbf{N}(\mathfrak{p})+1)(\mathbf{N}(\mathfrak{p})^2+1).$$
\end{lem}

\begin{prf} 
For $e_\mathfrak{p}=1$, the cardinality of the Lagrange variety over the finite field $\F_\mathfrak{p}=\mathcal{O}_F/\mathfrak{p}$ is given by $(\mathbf{N}(\mathfrak{p})+1)(\mathbf{N}(\mathfrak{p})^2+1)$.
Proceed by induction on $e_\mathfrak{p}$.
\end{prf}

{We have more to say about elements of $\mathcal{H}_0^2(\mathfrak{p}^{e_\mathfrak{p}})$ in Subsection~\ref{subsection:implementation}.}

\subsection{\bf Brandt matrices}  Let $\mathcal{F}=\left\{x_1,\,\ldots,\,x_h\right\}$ be a fundamental domain for the action of $\Gamma$ on $\mathcal{H}_0^2(N)$ and, for each $i$, let $\Gamma_i$ be the stabilizer of $x_i$. Then, every element in $M_k^B(N)$ is completely determined by its values on $\mathcal{F}$. Thus, there is an isomorphism of complex spaces
\begin{eqnarray*}
M_{k}^B(N)&\to&\bigoplus_{i=1}^h V_{k}^{\Gamma_i}\\
f&\mapsto& (f(x_i)),
\end{eqnarray*} where $V_{k}^{\Gamma_i}$ is the subspace of $\Gamma_i$-invariants in $V_{k}$.

For any $x,\,y\in \mathcal{H}_0^2(N)$, we let
$$
\begin{aligned}
\Theta_1(x,\,y,\,\mathfrak{p}) &:=\left\{u\in\Theta_1(\mathfrak{p})
\tq \exists \gamma\in\Gamma, \ ux=\gamma y \right\},\\
\Theta_2(x,\,y,\,\mathfrak{p}) &:=\left\{u\in\Theta_2(\mathfrak{p})
\tq \exists \gamma\in\Gamma, \ ux=\gamma y \right\}
\end{aligned}
$$

\begin{prop}\label{prop: brandt}
The actions of the Hecke operators $T_s(\mathfrak{p})$, $s=1,\,2$, are given by the Brandt matrices
$\mathcal{B}_s(\mathfrak{p})=(b_{ij}^s(\mathfrak{p}))$, where
\begin{eqnarray*}
b_{ji}^s(\mathfrak{p}):V_{k}^{\Gamma_j}&\to& V_{k}^{\Gamma_i}\\
v&\mapsto&v\cdot\left(\sum_{u\in\mathbf{\Theta}_s(x_i,\,x_j,\,\mathfrak{p})}\gamma_u^{-1}u\right).
\end{eqnarray*}
\end{prop}

\begin{prf}
The proof of Proposition~\ref{prop: brandt} follows the lines of \cite[\S 3]{dembele1}.
\end{prf}

\subsection{\bf Computing the group $G^B(\Z)$} It is enough to compute the subgroup $\Gamma$ consisting of the elements in $G^B(\Z)$ with similitude factor 1.  But it is easy to see that 
$$
\Gamma=\left\{\begin{pmatrix}u&0\\ 0&v\end{pmatrix} \ \Big\vert\ u,\,v\in\mathcal{O}_B^1\right\}\cup\left\{\begin{pmatrix}0&u\\ v&0\end{pmatrix} \ \Big\vert\ u,\,v\in\mathcal{O}_B^1\right\},
$$
where $\mathcal{O}_B^1$ is the group of norm 1 elements.

\subsection{\bf Computing the sets $\Theta_1(\mathfrak{p})$ and $\Theta_2(\mathfrak{p})$} Let us consider the quadratic form on the vector space $V=B^2$ given by
\begin{eqnarray*}
V&\to& F\\
(a,\,b)&\mapsto& |\!|(a,\,b)|\!|:=\mathbf{nr}(a)+\mathbf{nr}(b),
\end{eqnarray*} 
where $\mathbf{nr}$ is the reduced norm on $B$. This determines an inner form 
\begin{eqnarray*}
V\times V&\to& F\\
(u,\,v)&\mapsto& \langle u,\,v\rangle.
\end{eqnarray*} 
An element of  $\Theta_1(\mathfrak{p})$ (resp. $\Theta_2(\mathfrak{p})$) is a unitary matrix $\gamma\in \mathbf{M}_2(\mathcal{O}_B)$ with respect to this inner form such that the norm of each row is $\pi_\mathfrak{p}$ (resp. $\pi_\mathfrak{p}^2$ and the rank of the reduced matrix is $1$). So we first start by computing all the vectors $u=(a,\,b)\in \mathcal{O}_B^2$ such that $|\!|u|\!|=\pi_\mathfrak{p}$ (resp. $|\!|u|\!|=\pi_\mathfrak{p}^2$). And for each such vector $u$, we compute the vectors $v=(c,\, d)\in\mathcal{O}_B^2$ of the same norm such that $\langle u,\,v\rangle=0$. The corresponding matrix $\gamma=\begin{pmatrix}a&b\\ c&d\end{pmatrix}$ belongs to $\Theta_1(\mathfrak{p})$ (resp. $\Theta_2(\mathfrak{p})$) when its reduction mod $\mathfrak{p}$ has the appropriate rank. We list all these matrices up to equivalence and stop when we reach the right cardinality.

\subsection{\bf The implementation of the algorithm}\label{subsection:implementation}
The implementation of the algorithm is similar to that of \cite{dembele1}. However, it is important to note how we represent elements in $\mathcal{H}_0^2(N)$ so that we can retrieve them easily once stored. As in \cite{dembele1} we choose to
work with the product 
$$
\mathcal{H}_0^2(N)=\prod_{\mathfrak{p}\mid N}\mathcal{H}_0^2(\mathfrak{p}^{e_\mathfrak{p}}).
$$
Using Plucker's coordinates, we can view $\mathcal{H}_0^2(\mathfrak{p}^{e_\mathfrak{p}})$ as a closed subspace of
$\mathbf{P}^5(\mathcal{O}_{F_\mathfrak{p}}/\mathfrak{p}^{e_\mathfrak{p}})$. We then represent each element in $\mathcal{H}_0^2(\mathfrak{p}^{e_\mathfrak{p}})$ by choosing a point $x=(a_0:\cdots:a_5)=[u\wedge v]\in
\mathbf{P}^5(\mathcal{O}_{F_\mathfrak{p}}/\mathfrak{p}^{e_\mathfrak{p}})$ such that the submodule $M$ generated by $u$ and $v$ is a Lagrange submodule, and the first {\it invertible} coordinate is scaled to 1.

\begin{rem}\rm In \cite{lansky1}, Lansky and Pollack describe an algorithm which computes algebraic modular forms on the same inner form of $\mathbf{GSp}_4/\Q$ that we use. We would like to note that there are some differences between the two algorithms. Although \cite{lansky1} also uses the flag variety $\mathcal{H}_0^2(N)$ in order to determine the double coset space $G^B(\Q)\backslash G^B(\hat{\Q})/U_0(N)$, it later returns to the adelic setting in order to compute the Brandt matrices. In contrast, Theorem~\ref{thm1} and Proposition~\ref{prop: brandt} allow us to avoid that unnecessary step by describing the Hecke action on the flag variety $\mathcal{H}_0^2(N)$ directly. As a result, we get an algorithm that is more efficient.
\end{rem}

\section{\bf Numerical examples: $F=\Q(\sqrt{5})$ and $B=\left(\frac{-1,\,-1}{F}\right)$}\label{section: examples} 

In this section, we provide some numerical examples using the quadratic field $F=\Q(\sqrt{5})$. It is proven in K. Hashimoto and T. Ibukiyama \cite{hashimoto} that, for the Hamilton quaternion algebra $B$ over $F$, the class number of the principal genus of $G^B$ is one. We use our algorithm to compute all the systems of Hecke eigenvalues of Hilbert-Siegel cusp forms of weight 3 and level $N$ that are defined over real quadratic fields, where $N$ runs over all prime ideals of norm less than 50. We then determine which of the forms we obtained are possible lifts of Hilbert cusp forms by comparing the Hecke  eigenvalues for those primes.

\subsection{\bf Tables of Hilbert-Siegel cusp forms of parallel weight 3} In Table~\ref{hsmf} we list all the systems of eigenvalues of Hilbert-Siegel cusp forms of weight 3 and level $N$ that are defined over real quadratic fields, where $N$ runs over all prime ideals in $F$ of norm less than $50$. Here are the conventions we use in the tables.
\begin{enumerate}
\item For a quadratic field $K$ of discriminant $D$, we let $\omega_D$ be a generator of the ring of integers $\mathcal{O}_K$ of $K$.
\item The first row contains the level $N$, given in the format $(\text{Norm}(N), \alpha)$ for some generator $\alpha\in F$ of $N$, and the dimensions of the relevant spaces. 
\item The second row lists the Hecke operators that have been computed.
 \item For each eigenform $f$, the Hecke eigenvalues are given in a row, and the last entry of that row indicates if the form $f$ is a probable lift.
\item The levels and the eigenforms are both listed up to Galois conjugation.
\end{enumerate} 
For an eigenform $f$ and a given prime $\mathfrak{p}\nmid N$, let $a_1(\mathfrak{p},\,f)$ and $a_2(\mathfrak{p},\,f)$ be the eigenvalues of the Hecke operators $T_1(\mathfrak{p})$ and $T_2(\mathfrak{p})$, respectively. Then the Euler factor $L_\mathfrak{p}(f,\,s)$ is given (for example, in \cite[\S 3.4]{AS}) by
\[
L_\mathfrak{p}(f,s) =Q_\mathfrak{p}(q^{-s})^{-1},
\]
where
\begin{eqnarray*}
Q_\mathfrak{p}(x)&=&1- a_1(\mathfrak{p},f)x+b_1(\mathfrak{p}, f)x^2 - a_1(\mathfrak{p},f) q^{2k-3}x^3 + q^{4k-6}x^{4},\\
b_1(\mathfrak{p}, f)&=&a_1(\mathfrak{p},f)^2-a_2(\mathfrak{p}, f)-q^{2k-4}, \\
q &=& \mathbf{N}(\mathfrak{p}).
\end{eqnarray*}

\begin{table*}
\small
\begin{eqnarray*}\begin{array}{|crrrrrrc|}\hline
\multicolumn{8}{|l|}{N=(4,\,2):\quad\dim M_3^B(N)=2,\quad\dim S_3^B(N)=1}\\\hline
&T_1(2)&T_2(2)&T_1(\sqrt{5})&T_2(\sqrt{5})&T_1(3)&T_2(3)&\mbox{\rm Lift?}\\\hline\hline
f_{1}&-4&0&20&-36&140&580&yes
\\\hline
\multicolumn{8}{c}{}\\\hline
\multicolumn{8}{|l|}{N=(5,\,2+\omega_{5}):\quad\dim M_3^B(N)=2,\quad\dim S_3^B(N)=1}\\\hline
&T_1(2)&T_2(2)&T_1(\sqrt{5})&T_2(\sqrt{5})&T_1(3)&T_2(3)&\mbox{\rm Lift?}\\\hline\hline
f_{1}&20&15&-5&0&40&-420&yes
\\\hline
\multicolumn{8}{c}{}\\\hline
\multicolumn{8}{|l|}{N=(9,\,3):\quad\dim M_3^B(N)=3,\quad\dim S_3^B(N)=2}\\\hline
&T_1(2)&T_2(2)&T_1(\sqrt{5})&T_2(\sqrt{5})&T_1(3)&T_2(3)&\mbox{\rm Lift?}\\\hline\hline
f_{1}&25-3\omega_{41}&40-15\omega_{41}&30+6\omega_{41}&24+36\omega_{41}&-9&0&yes
\\\hline
\multicolumn{8}{c}{}\\\hline
\multicolumn{8}{|l|}{N=(11,\,3+\omega_{5}):\quad\dim M_3^B(N)=3,\quad\dim S_3^B(N)=2}\\\hline
&T_1(2)&T_2(2)&T_1(\sqrt{5})&T_2(\sqrt{5})&T_1(3)&T_2(3)&\mbox{\rm Lift?}\\\hline\hline
f_{1}&24&35&34&48&88&60&yes\\
f_{2}&-20&35&-10&4&0&60& no\\\hline
\multicolumn{8}{c}{}\\\hline
\multicolumn{8}{|l|}{N=(19,\,4+\omega_{5}):\quad\dim M_3^B(N)=5,\quad\dim S_3^B(N)=4}\\\hline
&T_1(2)&T_2(2)&T_1(\sqrt{5})&T_2(\sqrt{5})&T_1(3)&T_2(3)&\mbox{\rm Lift?}\\\hline\hline
f_{1}&4&11&-20&28&6&76&no\\
f_{2}&7&-50&15&-66&73&-90&yes\\
f_{3}&24+\omega_{161}&35+5\omega_{161}&36-\omega_{161}&60-6\omega_{161}&98-3\omega_{161}&160-30\omega_{161}&yes
\\\hline
\multicolumn{8}{c}{}\\\hline
\multicolumn{8}{|l|}{N=(29,\,5+\omega_{5}):\quad\dim M_3^B(N)=9,\quad\dim S_3^B(N)=8}\\\hline
&T_1(2)&T_2(2)&T_1(\sqrt{5})&T_2(\sqrt{5})&T_1(3)&T_2(3)&\mbox{\rm Lift?}\\\hline\hline
f_{1}&-4&11&10&20&30&60&no\\
f_{2}&8&-45&30&24&50&-320&yes\\
f_{3}&17&0&9&-102&86&40&yes
\\\hline
\multicolumn{8}{c}{}\\\hline
\multicolumn{8}{|l|}{N=(31,\,5+2\omega_{5}):\quad\dim M_3^B(N)=12,\quad\dim S_3^B(N)=11}\\\hline
&T_1(2)&T_2(2)&T_1(\sqrt{5})&T_2(\sqrt{5})&T_1(3)&T_2(3)&\mbox{\rm Lift?}\\\hline\hline
f_{1}&13&-20&20&-36&76&-60&yes
\\\hline
\multicolumn{8}{c}{}\\\hline
\multicolumn{8}{|l|}{N=(41,\,6+\omega_{5}):\quad\dim M_3^B(N)=19,\quad\dim S_3^B(N)=18}\\\hline
&T_1(2)&T_2(2)&T_1(\sqrt{5})&T_2(\sqrt{5})&T_1(3)&T_2(3)&\mbox{\rm Lift?}\\\hline\hline
f_{1}&10&20&-10&29&30&-20&no\\
f_{2}&-1&1 &5&14&-2&-56&no\\
f_{3}&27&50&40&84&124&420&yes\\
f_{4}&-12&19&30&65&0&0&no\\   
f_{5}&16-2\omega_{21}&-5-10\omega_{21} &21+4\omega_{21}&-30+24\omega_{21} &72-2\omega_{21}&-100-20\omega_{21}&yes\\
f_{6}&2-6\omega_{5}&11-2\omega_{5} &8+4\omega_{5}&11-4\omega_{5}&-12+54\omega_{5}&160+40\omega_{5}&no\\\hline
\multicolumn{8}{c}{}\\\hline
\multicolumn{8}{|l|}{N=(49,\,7):\quad\dim M_3^B(N)=26,\quad\dim S_3^B(N)=25}\\\hline
&T_1(2)&T_2(2)&T_1(\sqrt{5})&T_2(\sqrt{5})&T_1(3)&T_2(3)&\mbox{\rm Lift?}\\\hline\hline
f_{1}&5&-60&46&120&40&-420&yes\\
f_{2}&4+4\omega_{65}&32+3\omega_{65}&12-4\omega_{65}&44-4\omega_{65}&-6-12\omega_{65}&145+8\omega_{65}&no\\ \hline
\end{array}\end{eqnarray*}
\caption{\bf Hilbert-Siegel eigenforms of weight 3}
\label{hsmf}
\end{table*}

\subsection{\bf Tables of Hilbert cusp forms of parellel weight 4}
In Table~\ref{hmf}, we list all the Hilbert cusp forms of parallel weight $4$ and level $N$ that are defined over real quadratic fields, with $N$ running over all prime ideals of norm less than 50. (They are computed by using the algorithm in \cite{dembele1}). We use this data in order to determine the forms in Table~\ref{hsmf} that are possible lifts from $\mathbf{GL}_2$. 

\begin{table*}
\small
\begin{eqnarray*}\begin{array}{|c|c|c|c|c|c|c|c|c|c|}\hline 
\multicolumn{2}{|c|}{N}&\multicolumn{1}{|c|}{(4,2)}&\multicolumn{1}{|c|}{(5,2+\omega_{5})}&\multicolumn{1}{|c|}{(9,3)}&\multicolumn{1}{|c|}{(11,3+\omega_{5})}\\\hline
\mathbf{N}(\mathfrak{p})&\mathfrak{p}&a(\mathfrak{p},\, f_{1})&a(\mathfrak{p},\, f_{1})&a(\mathfrak{p},\, f_{1})&a(\mathfrak{p},\, f_{1})\\\hline\hline
4&2&-4&0&5-3\omega_{41}&4\\ 
5&2+\omega_{5}&-10&-5&6\omega_{41}&4\\ 
9&3&50&-50&-9&-2\\ 
11&3+2\omega_{5}&-28&32&-18-6\omega_{41}&-10\\ 
11&3+\omega_{5}&-28&32&-18-6\omega_{41}&-11\\ 
19&4+3\omega_{5}&60&100&-40+24\omega_{41}&-94\\ 
19&4+\omega_{5}&60&100&-40+24\omega_{41}&28\\ 
\hline
\end{array}\end{eqnarray*}
\begin{eqnarray*}\begin{array}{|c|c|c|c|c|c|c|c|c|c|}\hline 
\multicolumn{2}{|c|}{N}&\multicolumn{2}{|c|}{(19,4+\omega_{5})}&\multicolumn{2}{|c|}{(29,5+\omega_{5})}\\\hline
\mathbf{N}(\mathfrak{p})&\mathfrak{p}&a(\mathfrak{p},\, f_{1})&a(\mathfrak{p},\, f_{2})&a(\mathfrak{p},\, f_{1})&a(\mathfrak{p},\, f_{2})\\\hline\hline
4&2&-13&5-\omega_{161}&-12&-3\\ 
5&2+\omega_{5}&-15&5+\omega_{161}&0&-21\\ 
9&3&-17&5+3\omega_{161}&-40&-4\\ 
11&3+2\omega_{5}&-6&2+8\omega_{161}&-68&37\\ 
11&3+\omega_{5}&33&7-7\omega_{161}&30&-66\\ 
19&4+3\omega_{5}&-139&-15-9\omega_{161}&-28&-40\\ 
19&4+\omega_{5}&19&-19&84&-9\\ 
\hline
\end{array}\end{eqnarray*}
\begin{eqnarray*}\begin{array}{|c|c|c|c|c|c|c|c|}\hline 
\multicolumn{2}{|c|}{N}&\multicolumn{1}{|c|}{(31,5+2\omega_{5})}&\multicolumn{2}{|c|}{(41,6+\omega_{5})}\\\hline
\mathbf{N}(\mathfrak{p})&\mathfrak{p}&a(\mathfrak{p},\, f_{1})&a(\mathfrak{p},\, f_{1})&a(\mathfrak{p},\, f_{2})\\\hline\hline
4&2&-7&7&-4-2\omega_{21}\\ 
5&2+\omega_{5}&-10&10&-9+4\omega_{21}\\ 
9&3&-14&34&-18-2\omega_{21}\\ 
11&3+2\omega_{5}&-20&-60&-19\\ 
11&3+\omega_{5}&-28&-2&-24-4\omega_{21}\\ 
19&4+3\omega_{5}&-12&74&4-50\omega_{21}\\ 
19&4+\omega_{5}&28&16&-29+44\omega_{21}\\ 
\hline
\end{array}\end{eqnarray*}
\begin{eqnarray*}\begin{array}{|c|c|c|c|c|c|}\hline 
\multicolumn{2}{|c|}{N}&\multicolumn{2}{|c|}{(49,7)}\\\hline
\mathbf{N}(\mathfrak{p})&\mathfrak{p}&a(\mathfrak{p},\, f_{1})&a(\mathfrak{p},\, f_{2})\\\hline\hline
4&2&-15&-2\\ 
5&2+\omega_{5}&16&-10\\ 
9&3&-50&-11\\ 
11&3+2\omega_{5}&-8&-7-28\omega_{13}\\ 
11&3+\omega_{5}&-8&-35+28\omega_{13}\\ 
19&4+3\omega_{5}&-110&-26+14\omega_{13}\\ 
19&4+\omega_{5}&-110&-12-14\omega_{13}\\ 
\hline
\end{array}\end{eqnarray*}
\caption{\bf Hilbert eigenforms of weight 4}
\label{hmf}
\vspace{-1cm}
\end{table*}

\subsection{\bf Lifts} 
There are two types of lifts from $\mathbf{GL}_2$ to $\mathbf{GSp}_4$. The first one corresponds to the homomorphism of $L$-groups determined by the long root embedding into $\mathbf{GSp}_4$, and the second one by the short root embedding. (See \cite{lansky1} for more details). Let $f$ be a Hilbert cusp form of parallel weight $k$ and level $N$ with Hecke eigenvalues $a(\mathfrak{p},\,f)$, where $\mathfrak{p}$ is a prime not dividing $N$. Let $\phi$ be the lift of $f$ to $\mathbf{GSp}_4$ via the long root, and $\psi$ the one via the short root. Then the Hecke eigenvalues of $\phi$ are given by
\[
\begin{aligned}
a_1(\mathfrak{p},\,\phi)&=a(\mathfrak{p},\,f)\ \mathbf{N}(\mathfrak{p})^{\frac{4-k}{2}}+\mathbf{N}(\mathfrak{p})^2+\mathbf{N}(\mathfrak{p})\\
a_2(\mathfrak{p},\,\phi)&=a(\mathfrak{p},\,f)\ \mathbf{N}(\mathfrak{p})^{\frac{4-k}{2}}(\mathbf{N}(\mathfrak{p})+1)+\mathbf{N}(\mathfrak{p})^2-1,
\end{aligned}
\]
and the Hecke eigenvalues of $\psi$ are given by
\[
\begin{aligned}
a_1(\mathfrak{p},\,\psi)&=a(\mathfrak{p},\,f)^3\ \mathbf{N}(\mathfrak{p})^{\frac{6-3k}{2}}-2\ a(\mathfrak{p},\,f)\ \mathbf{N}(\mathfrak{p})^{\frac{4-k}{2}}\\
a_2(\mathfrak{p},\,\psi)&=a(\mathfrak{p},\,f)^4\ \mathbf{N}(\mathfrak{p})^{4-2k}-3\ a(\mathfrak{p},\,f)^2\ \mathbf{N}(\mathfrak{p})^{3-k}+\mathbf{N}(\mathfrak{p})^2-1.
\end{aligned}
\]
The second lift $\psi$ is the so-called symmetric cube lifting.

\begin{rem}\rm So far, our algorithm has been implemented only for congruence subgroups of Siegel type. We intend to improve the implementation in the near future so as to include more additional level structures such as the Klingen type. Indeed, Ramakrishnan and Shahidi \cite{dinakar} recently showed the existence of symmetric cube lifts for non-CM elliptic curves $E/\Q$ to $\mathbf{GSp}_4/\Q$. And their result should hold for other totally real number fields, with the level structures of the lifts being of Klingen type. Unfortunately, those lifts cannot be seen in our current tables. For example, there are modular elliptic curves over $\Q(\sqrt{5})$ whose conductors have norm 31, 41 and 49, but the corresponding symmetric cubic lifts do not appear in Table~\ref{hsmf}. We would like to remedy that in our next implementation.
\end{rem}


\end{document}